\documentclass[12pt]{article}
\usepackage{fullpage}
\usepackage{amssymb,amsmath,amsfonts, amscd}
\usepackage{array}
\usepackage{amsmath}
\usepackage{amssymb}

\begin{document}

\newtheorem{theorem}{Theorem} [section]
\newtheorem{conjecture}[theorem]{Conjecture}
\newtheorem{tha}{Theorem}
\renewcommand{\thetha}{\Alph{tha}}
\newtheorem{corollary}[theorem]{Corollary}
\newtheorem{lemma}[theorem]{Lemma}
\newtheorem{proposition}[theorem]{Proposition}
\newtheorem{construction}[theorem]{Construction}
\newtheorem{question}[theorem]{Question}
\newtheorem{definition}[theorem]{Definition}
\newtheorem{observation}{Observation}
\newtheorem{remark}[theorem]{Remark}
\newtheorem{fact}[theorem]{Fact}
\setlength{\textwidth}{17cm}
\setlength{\oddsidemargin}{-0.1 in}
\setlength{\evensidemargin}{-0.1 in}
\setlength{\topmargin}{0.0 in}

\def \df {\noindent {\bf Definition. }}

\newcommand{\dist}{{\rm dist}}
\newcommand{\Dist}{{\rm Dist}}
\newcommand{\chib}{\chi_B}
\newcommand{\Forb}{{\rm Forb}}

\newcommand{\qedbox}{$\blacksquare$ \newline}
\newenvironment{proof}%
{%
\noindent{\it Proof.} } {
\hfill\qedbox }

\newcommand{\proofend}{\hfill\qedbox}

\def\qed{\hskip 1.3em\hfill\rule{6pt}{6pt} \vskip 20pt}

\linespread{1.0}
\input epsf
\def\epsfsize#1#2{1.0#1\relax}
\def\O{\text{O}}
\def\o{\text{o}}
\def\ex{\text{ex}}
\def\Z{\mathbb Z}

\def\fl#1{\lfloor #1 \rfloor}
\def\ce#1{\lceil #1 \rceil}

\def \cH{{\cal H }}
\def \cF{{\cal F}}
\def \cG{{\cal G}}
\def \cQ{{\cal Q}}
\def \cA{{\cal A}}
\def \cD{{\cal D}}

\def \f {{\cal F }}
\def \A {{\cal A}}
\def \D {{\cal D}}
\def \fn2 {{\lfloor n/2 \rfloor}}
\def \cn2 {{\lceil  n/2 \rceil}}

\newcommand{\bin}[2]{{#1\choose #2}}
\newcommand{\comp}{\overline}
\newcommand{\exval}{{\rm ex}}

\newcommand{\dotcup}{\,\dot\cup\,}

\title{On the editing distance of graphs}
\author{Maria Axenovich\thanks{Department of Mathematics, Iowa
State University, Ames, IA 50011, {\tt axenovic@math.iastate.edu}}
\and Andr\'e K\'ezdy\thanks{Department of Mathematics, University
of Louisville, Louisville, KY 40292, {\tt kezdy@louisville.edu}}
\and Ryan Martin\thanks{Department of Mathematics, Iowa State
University, Ames, IA 50011, {\tt
rymartin@iastate.edu}}~\thanks{Corresponding author.}}

\date{\today}
\maketitle

\begin{abstract}
An edge-operation on a graph $G$ is defined to be either the
deletion of an existing edge or the addition of a nonexisting
edge. Given a family of graphs $\cG$, the editing distance from
$G$ to $\cG$ is the smallest number of edge-operations needed to
modify $G$ into a graph from $\cG$. In this paper, we fix a graph
$H$ and consider $\Forb(n,H)$, the set of all graphs on $n$
vertices  that have no induced copy of $H$. We provide bounds for
the maximum over all $n$-vertex graphs $G$ of the editing distance
from $G$ to $\Forb(n,H)$, using an invariant we call the {\it
binary chromatic number} of the graph $H$. We give asymptotically
tight bounds for that distance when $H$ is self-complementary and
exact results for several small graphs $H$.
\end{abstract}

\section{Introduction}
\label{intro}

The investigation of graphs not containing subgraphs with given
properties is a classical problem. For example, determining the
maximum number of edges in a graph with no copy of a fixed
subgraph $H$ has been studied intensively for the last 70 years
\cite{F,S}. Very often, though, the desired task is not to
determine the extremal graph without  a given fixed subgraph, but
rather to start with an arbitrary graph and modify it in a small
number of steps such that the resulting graph does not contain a
forbidden subgraph.

The problem of modifying the given graph such that the resulting
graph satisfies some {\it global properties} has been addressed by
Erd\H{o}s et al. \cite{EGR,EFPS, EGS}. They investigated the number of
edge deletions sufficient to transform an arbitrary triangle-free
graph into a bipartite graph, as well as the smallest number of
edge additions sufficient to decrease diameter.

In this paper, we investigate the problem of transforming a given
graph into a new graph having the {\it local property} of avoiding
a fixed induced subgraph. Starting with an arbitrary graph $G$, we
would like to calculate the minimum number of edges needed to be
added to or deleted from $G$ to obtain a graph not containing a
fixed induced subgraph. Formally, an {\bf edge-operation} is
defined to be either the deletion of an existing edge or the
addition of a nonexisting edge.  Let $\Dist(G,H)$ denote the
minimum number of edge-operations needed to transform the graph
$G$ into a graph isomorphic to $H$.
In other words, it can be described as a symmetric difference as follows:
 $$\Dist(G,H)= \min \{|E(G)\Delta E(H')|~:~H'\cong H\}.$$
 Clearly this parameter is
defined if and only if $G$ and $H$ have the same number of
vertices.  If ${\cal H}$ is a class of graphs on $n$ vertices, we
define $\Dist (G, {\cal H})= \min\{ \Dist(G,H):\, H\in {\cal H}
\}$ for a graph $G$ on $n$ vertices. Finally, $\Dist(n, {\cal
H})=\max\{\Dist (G, {\cal H}): \; |V(G)|= n\}$. We call the metric
$\Dist(G,H)$ the {\bf editing distance} since the operations
performed can be considered as editing the edge set of a graph.
Our interest here is the class ${\cal H}$
of graphs on $n$ vertices containing no copies of a given fixed
graph $H$ as an induced subgraph. We denote this class by
$\Forb(n,H)$ (or simply by $\Forb(H)$ when it is clear from the
context). Similarly $\Forb'(H)$ is the family of all graphs on $n$
vertices with no subgraph isomorphic to $H$.

This graph editing problem has numerous applications in computer
science and bioinformatics. For example,  consider a metabolic
network and identify genes with vertices of a graph and pairs of
interacting genes with edges of the graph. It is a fundamental
question in biology (from evolutionary and practical points of
view) to find how many edge-changes in such a graph must be
performed to avoid an induced subgraph corresponding to a certain
metabolic process. Another example involves consensus trees. It is
known that two consensus trees are comparable if there is no
induced path on five vertices in a  corresponding bipartite
graph~\cite{SAN1998,CEFS,DF}. In particular, finding the smallest
number of edge-changes in such a graph will determine the distance
between these trees.

On the other hand, the editing problem of graphs corresponds to
determining the distance between $\{0,1\}$-matrices. If $A$ and $B$
are the adjacency matrices of graphs  $G_1$ and $G_2$ respectively,
then $\Dist(G_1,G_2)$ corresponds to the number of positions where
$A$ and $B$ differ, i.e., to the  Hamming distance between $A$ and
$B$. Thus finding editing distance between classes of graphs
provides the Hamming distance between classes of symmetric
matrices with the same diagonal entries. Moreover, when the graph
editing problem is restricted to bipartite graphs in which edge
additions and deletions are limited to edges between partite sets,
it corresponds to the problem of determining the distance between
the sets of arbitrary $\{0,1\}$-matrices.

We define the distance $\Dist'(n, \Forb'(H))$ to be analogous to
$\Dist(n,\Forb(H))$, but in this case, only permit edge-deletions.
 This quantity will always be equal to $\Dist(K_n, \Forb'(H))$, i.e.,
it is the minimum number of edges in the complement of an $H$-free
graph. If $\exval(n,H)$ is the maximum number of edges in a graph
on $n$ vertices with no subgraph isomorphic to $H$, then
\begin{equation}
   \Dist(K_n, \Forb'(H))=\binom{n}{2}-\exval(n, H) .
   \label{turanex}
\end{equation}
The asymptotic behavior of $\exval(n,H)$ is provided by the
following theorem, which was generalized by Erd\H{o}s and
Simonovits \cite{ES2}.
\begin{theorem}[Erd\H{o}s, Stone~\cite{ES}]
$\exval(n,H)=\left(1-\frac{1}{\chi(H)-1}+o(1)\right)
\binom{n}{2}$.
\end{theorem}
In particular, the distance $\Dist(n, \Forb'(H))$ is
asymptotically determined by the chromatic number of a graph $H$.

Clearly, when a forbidden graph is complete or empty, finding the
editing distance becomes a trivial task immediately reduced to
Tur\'an's theorem \cite{T}. On the other hand, perhaps the most
interesting case is when the forbidden induced subgraph is
self-complementary, i.e., when both operations of edge-deletions
and edge-additions carry ``an equal power''. In this case, we
derive  asymptotically tight estimates for $\Dist(n, \Forb(H))$.
We also give general bounds for other graphs. Our main tool in
providing the lower bounds is Szemer\'edi's Regularity Lemma which
allows us to express the bounds in terms of an invariant which we
designate the {\bf binary chromatic number}.  In defining it, we
use the term {\bf coclique} in place of the term {\bf independent
set}.
\begin{definition}
  The {\bf binary chromatic number} of a graph $G$, $\chib(G)$
  is the least integer $k+1$ such that, for all
  $c\in\{0,\ldots,k+1\}$, there exists a partition of $V(G)$
  into $c$ cliques and $k+1-c$ cocliques.
\end{definition}

This invariant was first introduced by Pr\"omel and Steger (called
$\tau$ in \cite{PS}) to express, asymptotically, the number of
$n$-vertex graphs which fail to have an induced copy of some
small, fixed graph $H$.  In \cite{BT,B}, $\tau$ was generalized as
the so-called colouring number of a hereditary property, $\cal P$.
In particular, when ${\cal P}=\Forb(H)$, the colouring number of
$\cal P$ is exactly $\chib(H)-1$. It should be noted that $\chib$
is not the cochromatic number (see \cite{LS}) even though the
definitions may seem to be similar at first glance.  We use the
term ``binary chromatic number'' in order to emphasize the close
connection to the chromatic number and the complementary r\^ole
that cliques and cocliques play.

The following comprise the main results of this paper.

\begin{theorem}
If $H$ is a graph with binary chromatic number $k+1$, then
$$ \Dist(n,\Forb(H))>(1-o(1))\frac{n^2}{4k}. $$
\label{tLB}
\end{theorem}

If $k=\chib(G)-1$, then let $c_{\min}$ be the least $c$ so that
$G$ cannot be partitioned into $c$ cliques and $k-c$ cocliques.
Let $c_{\max}$ be the greatest such number. We now have an upper
bound that can be expressed in terms of the binary chromatic
number of $H$ and corresponding $c_{\min}$ and $c_{\max}$.
\begin{theorem}
   Let $H$ be a graph with binary chromatic number $k+1$
   and $c_{\min}$ and $c_{\max}$ be defined as
   above.  If $c_{\min}\leq k/2\leq c_{\max}$, then
   \begin{equation}
      \Dist(n,\Forb(H)) \leq \frac{1}{2k}\bin{n}{2} . \label{UB1}
   \end{equation}
   Otherwise, let $c_0$ be the one of $\{c_{\max},c_{\min}\}$ that
   is closest to $k/2$.  Then
   \begin{equation}
      \Dist(n,\Forb(H)) \leq \left(\frac{1}{1+
      2\sqrt{\frac{c_0}{k}\left(1-\frac{c_0}{k}\right)}}
      \right)\frac{1}{k}\bin{n}{2}\leq
      \frac{1}{k}\bin{n}{2} . \label{UB2}
   \end{equation}
   \label{tUB}
\end{theorem}

\begin{corollary}
   If $H$ is a self-complementary graph with the property that
   $\chib(H)=k+1$, then
   $$ \Dist(n,\Forb(H))=(1+o(1))\frac{n^2}{4k} . $$
   \label{maincor}
\end{corollary}

In  Section~\ref{defs} we give preliminary
definitions and results that aid in the proof of the theorems.
Section~\ref{proofs} contains proofs of the main theorems. We
investigate the properties of the binary chromatic number in
Section~\ref{chromnum}. Finally, Section~\ref{exact} gives several
exact results.

\section{Definitions and preliminary results}
\label{defs}

We denote by $K_n$, $E_n$, $C_n$, and $P_n$ a complete graph, an empty graph, a
cycle, and a path on $n$ vertices, respectively.  We also define
$K_p^q$ to be a complete $p$-partite graph with each partite set
of cardinality $q$.  We use $\comp{G}$ to denote the complement of
$G$.  For the other definitions, we refer the reader to \cite{W}.

\begin{definition} For a graph $G=(V,E)$ and two disjoint subsets $A$ and $B$ of vertices,
the {\bf density} of a pair $(A,B)$ is denoted $d(A,B)$ and is given by the
formula
$$ d(A,B)=\frac{e(A,B)}{|A||B|} $$
where $e(A,B)$ is the number of edges of $G$ with one end-point in $A$ and another in $B$.
\end{definition}

\begin{definition}
   For a graph $G=(V,E)$ and two disjoint subsets $A$ and $B$ of vertices, a pair $(A,B)$ is  {\bf $\epsilon$-regular} if
   $$ X\subset A, Y\subset B, |X|>\epsilon|A|, |Y|>\epsilon |B| $$
   imply
   $$ |d(X,Y)-d(A,B)|<\epsilon ; $$
   otherwise, $(A,B)$ is {\bf $\epsilon$-irregular}.
\end{definition}

\noindent The proof of Theorem~\ref{ForcedH} makes use of the
Regularity Lemma (see~\cite{KS} and~\cite{KSSSz}).

\begin{lemma}[Regularity Lemma~\cite{Sz}]\label{reglem}
For every positive
$\epsilon$ and positive integer $m$, there are positive integers
$M=M(\epsilon,m)$ and $N=N(\epsilon,m)$ with the following
property:  For every graph $G$ with at least $N$ vertices there is
a partition of the vertex set into $\ell+1$ classes (clusters)
$$ V=V_0\dotcup V_1\dotcup V_2\dotcup\cdots\dotcup V_{\ell} $$ such
that
\begin{enumerate}
   \item $m\leq\ell\leq M$, \label{reglem1}
   \item $|V_1|=|V_2|=\cdots=|V_{\ell}|$, \label{reglem2}
   \item $|V_0|<\epsilon n$, \label{reglem3}
   \item at most $\epsilon\ell^2$ of the pairs $(V_i,V_j)$ are
   $\epsilon$-irregular. \label{reglem4}
\end{enumerate}
\end{lemma}

\noindent
We give a name to the partition given by Lemma~\ref{reglem}:
\begin{definition}
   An {\bf $(m,\epsilon,\ell)$-equipartition} of a vertex
   set $V$ is a partition
   $V=\linebreak V_0\dotcup V_1\dotcup\cdots\dotcup V_{\ell}$ such
   that the Regularity Lemma's conditions (\ref{reglem1}),
   (\ref{reglem2}) and (\ref{reglem3}) are satisfied (with
   $M=M(\epsilon,m)$ as defined  by the Lemma).
\end{definition}

In order to state our lower bound, we need to generalize the idea
of an $\epsilon$-regular pair.
\begin{definition}
   An {\bf $\epsilon$-regular $r$-tuple} is an $r$-partite graph
   with partite sets $V_1,\ldots,V_r$ such that $|V_i|=|V_j|$
   and $(V_i,V_j)$ is an $\epsilon$-regular pair for all $i,j$,
   with $1\leq i<j\leq r$.

   We say that an $\epsilon$-regular $r$-tuple is {\bf of size
   $rL$} if $|V_1|=\cdots=|V_r|=L$.  For $0<\delta<1/2$, an
   $\epsilon$-regular $r$-tuple has {\bf $\delta$-bounded
   density} if $d(V_i,V_j)\in(\delta,1-\delta)$ whenever
   $1\leq i<j\leq r$.

   For convenience, we define an {\bf $(\epsilon, r, L,
   \delta)$-configuration} to be an $\epsilon$-regular $r$-tuple of
   size $rL$ that has $\delta$-bounded density.
   \label{delbounded}
\end{definition}

The following Theorem is our major tool in proving the main result. We prove it in Section \ref{proofs}.

\begin{theorem}
   Let $r$ be a positive integer and $\delta$ and $\epsilon$ be
   real numbers, $0<\delta<1$ and $\epsilon>0$, such that
   $\epsilon<\delta/(16r-16)$. There is a graph $G$ on $n$
   vertices and a constant $M(\epsilon)$ such that if the
   number of edge-deletions and edge-additions performed
   on $G$ is less than
   $\frac{n^2}{4(r-1)}(1-3\delta)(1-\epsilon)^2$
   then the resulting graph contains an
   $(\epsilon,r,L,\delta)$-configuration with
   $L\geq n(1-\epsilon)/M(\epsilon)$.
   \label{EditOperations}
\end{theorem}

The following theorem is essentially Lemma 3.5 in \cite{PS}. We shall use it to prove
Theorem \ref{tLB}.

\begin{theorem}[Pr\"omel, Steger~\cite{PS}]
   Let $H$ be a fixed graph with binary chromatic number $r$ and
   $\delta$ be a real number with $0<\delta<1/2$.  There
   exists an $\epsilon_0=\epsilon_0(H,\delta,r)>0$ such that for all
   $\epsilon$, where $0<\epsilon\leq\epsilon_0$, there exists an
   $n_0=n_0(H,\delta,\epsilon,r)$ such that every graph $G=(V,E)$ on
   $n\geq n_0$ vertices has the following property: Let
   $V=V_0\dotcup V_1\dotcup\cdots\dotcup V_{\ell}$ be an
   $(r,\epsilon,\ell)$-equipartition with
   $|V_1|=\cdots=|V_{\ell}|=L$ for $L\geq n/M(\epsilon,r)$
   where $M(\epsilon,r)$ is the constant given in Lemma
   \ref{reglem} for $\epsilon$ and $r$, such that
   $(V_1,\ldots,V_r)$ forms an
   $(\epsilon,r,L,\delta)$-configuration.  Then, the
   subgraph induced by $\bigcup_{i=1}^r V_i$ contains the
   graph $H$ as an induced subgraph.
   \label{ForcedH}
\end{theorem}
For functions $f=f(n)$ and $g=g(n)$, let $f=\omega(g)$ be the usual asymptotic notation
denoting that $g/f\rightarrow 0$ as $n\rightarrow\infty$.

\begin{lemma}
  For a positive real number  $\epsilon$ and positive integer $m$,
  let $\ell$ have the property that $m\leq\ell\leq M(\epsilon,m)$,
  where $M(\epsilon,m)$ is the constant given in the Regularity
   Lemma (Lemma \ref{reglem}). Let $f=f(n)=\omega(n^{-1/2})$.
Then, for $n$ large enough, there is
a graph $G$ on $n$ vertices so that for any
$(m,\epsilon,\ell)$-equipartition, all pairs of clusters
$(V_i,V_j)$, $1\leq i<j\leq \ell$, are $\epsilon$-regular with
density in the interval $(1/2- f, 1/2+f)$. \label{cRANDOM}
\end{lemma}

The proof of this Lemma is  a  routine calculation
that  we include in the Appendix,  for completeness.

\section{The proofs}
\label{proofs}

In order to prove Theorem \ref{tLB}, we prove Theorem \ref{EditOperations}, which basically
asserts that a graph described in Lemma \ref{cRANDOM}  requires many editing operations
to eliminate all induced copies of $H$.


\subsection{Proof of Theorem \ref{EditOperations}}

Fix $\delta$ such that $0<\delta<1/2$.  Let $G$ be a graph on $n$
vertices as described in Lemma~\ref{cRANDOM} with $\epsilon$,
$r$ and $L$ as given in that Lemma. Let $G'$ be a graph
with no $(\epsilon, r, L, \delta)$-configuration having least
distance from $G$.

Apply the Regularity Lemma to $G'$ with parameters $\epsilon$ and
$m=\epsilon^{-1}$, to get $\ell+1$ clusters (we have $\ell\geq
\epsilon^{-1}$) $V_0,V_1,\ldots,V_{\ell}$ such that
$|V_1|=\cdots=|V_{\ell}|=L$. Furthermore, all but $\epsilon\ell^2$
pairs $(V_i,V_j)$, $1\leq i<j\leq\ell$ are $\epsilon$-regular.
Recalling Definition \ref{delbounded}, we say that an
$\epsilon$-regular pair is {\bf $\delta$-bounded} if its density
is at least $\delta$ and at most $1-\delta$; otherwise, it is {\bf
$\delta$-unbounded}.

Since $G'$ does not have an $(\epsilon, r, L,
\delta)$-configuration, it is not possible to have a set of $r$
clusters such that between any two clusters,
there is a $\delta$-bounded $\epsilon$-regular pair. Thus,
according to Tur\'an's theorem, the number of pairs of clusters,
$(V_i,V_j)$, that induce either a $\delta$-unbounded
$\epsilon$-regular pair or an $\epsilon$-irregular pair is at
least
$$ (r-1)\frac{\frac{\ell}{r-1}\left(\frac{\ell}{r-1}-1\right)}{2}
   =\frac{\ell(\ell-r+1)}{2(r-1)} . $$

The number of $\epsilon$-irregular pairs, $(V_i,V_j)$, is at
most $\epsilon\ell^2$, thus the number of $\delta$-unbounded
$\epsilon$-regular pairs is at least
$$ \ell^2\left(\frac{1}{2(r-1)}-\epsilon\right)-\frac{\ell}{2} .
$$

Because $G$ came from Lemma~\ref{cRANDOM}, if some pair
$(V_i,V_j)$ were $\delta$-unbounded in $G'$, then at least
$(1/2-\delta-o(1))L^2$ edges had to be either added or deleted between $V_i$ and $V_j$ in order to get $G'$ from $G$.
Hence, the total number of edges that had to be changed is at
least
$$ \ell^2L^2
   \left(\frac{1}{2(r-1)}-\epsilon-\frac{1}{2\ell}\right)
   \left(\frac{1}{2}-\delta-o(1)\right)\geq
   \frac{\ell^2L^2}{4(r-1)}(1-3\delta) . $$
The inequality is valid as long as
$8(r-1)(\epsilon+\ell^{-1}/2)<\delta$.

Since $\ell L\geq n(1-\epsilon)$, the total number of edges that
have to be  altered  in order to obtain $G'$ from $G$ is at least
$$ \frac{n^2}{4(r-1)}(1-3\delta)(1-\epsilon)^2 . $$
\proofend

\subsection{Proof of Theorem~\ref{tLB}}

Choose a $\delta$ arbitrarily small, and let $G$ be the graph
guaranteed by Theorem~\ref{EditOperations}.  If fewer than
$\frac{n^2}{4k}(1-3\delta)(1-\epsilon)^2$ edge-operations are
performed on $G$ to obtain a graph $G'$, then there are disjoint  vertex sets
$V_1,\ldots,V_{k+1}$ in $G'$ that satisfy the conditions of
Theorem~\ref{ForcedH}. Theorem~\ref{ForcedH} then implies that
$G'$ contains an induced $H$. Thus the editing distance $\Dist(G, \Forb(H))$ is at least
$\frac{n^2}{4k}(1-3\delta)(1-\epsilon)^2$.
\proofend


\subsection{Proof of Theorem~\ref{tUB}}

Lemma~\ref{lUB} emphasizes the importance of $c$ in the definition
of $\chib$.
\begin{lemma}
   Let $H$ be a graph with binary chromatic number $k+1$ and
   $c$ be an integer, $0\leq c\leq k$, so that $H$ cannot
   be covered by
   exactly $c$ cliques and  $k-c$ independent sets.
   Let $G$ be a graph with density $d=e(G)/\bin{n}{2}$.  As long
   as it is not the case that both $d=0$ and $c=k$,
   or both $d=1$ and $c=0$,
   \begin{equation}
      \Dist(G,\Forb(H))\leq
      \frac{d(1-d)}{dc+(1-d)(k-c)}\bin{n}{2} ;
      \label{eqUB}
   \end{equation}
   otherwise,
   $\Dist(G,\Forb(H))\leq\frac{1}{k}\bin{n}{2}$.
   \label{lUB}
\end{lemma}

\begin{proof}
In order to prove the statement of the Lemma, we provide a
probabilistic algorithm adding and deleting some edges of $G$ such
that the resulting graph has no induced copy of $H$. We begin by
assigning colors independently to the vertices of $G$:
$1,\ldots,c$ each with probability $p$ and $c+1,\ldots,k$ each
with probability $q$.  Call such a coloring $g$. If $g(x)=g(y) \in
\{1,\ldots, c\}$ and $xy\notin E(G)$, then add an edge $xy$ to
$E(G)$. If $g(x)=g(y)\in \{c+1,\ldots, k\}$, and $xy \in E(G)$,
then delete $xy$ from $E(G)$. As a result, we obtain a graph $G'$
with the vertex set partitioned into $k$ subsets.  The first $c$
of these subsets induce cliques and the others induce cocliques.
Since the vertices of $H$ can not be partitioned into $c$ cliques
and $k-c$ cocliques, $H$ is not an induced subgraph of $G'$.

The expected number of changes is
$$ f(p,q)= \left(\bin{n}{2}-e(G)\right)c p^2+e(G)(k-c)q^2=
   \left((1-d)c p^2+d(k-c) q^2\right)\bin{n}{2} . $$
We also have the restriction
\begin{equation}
   c p+(k-c) q=1 . \label{g}
\end{equation}

As long as we do not have the case that both $d=0$ and $c=k$ or
the case that both $d=1$ and $c=0$, the method of Lagrange
multipliers gives that the minimum of $f(p,q)$ restricted to
(\ref{g}) occurs when $p=d/(dc+(1-d)(k-c))$ and
$q=(1-d)/(dc+(1-d)(k-c))$ and is equal to
$$ \frac{d(1-d)}{dc+(1-d)(k-c)}\bin{n}{2}. $$

Since this is the expected number of changes, there exists a
partition of the vertices of $G$ such that the above procedure
requires at most  $\frac{d(1-d)}{dc+(1-d)(k-c)}\bin{n}{2}$ changes to
make the graph $H$-free.

If both $d=0$ and $c=k$, then perform the above procedure, but fix
$p=1$.  If both $d=1$ and $c=0$, then perform the above procedure,
but fix $q=1$.  In both cases, the expected number of changes to
be performed is $\frac{1}{k}\bin{n}{2}$.
\end{proof}

In order to prove inequality (\ref{UB1}) of Theorem~\ref{tUB}, we
use Lemma \ref{lUB} and
find conditions when
\begin{equation}
   \frac{d(1-d)}{dc+(1-d)(k-c)}\leq\frac{1}{2k} . \label{ck}
\end{equation}

If $c\leq k/2$, then (\ref{ck}) holds when $d\in [0,1/2]\cup
[1-c/k,1]$.  If $c\geq k/2$, then (\ref{ck}) holds when $d\in
[0,1-c/k]\cup [1/2,1]$.  Consider a graph $G$ of density $d$ and
an $H$ for which $c_{\min}\leq k/2\leq c_{\max}$. If $d\leq 1/2$,
then choose $c_{\min}$; otherwise, choose $c_{\max}$.
 As a result,
$\Dist(G,\Forb(H))\leq\frac{1}{2k}\bin{n}{2}$.

In order to prove inequality (\ref{UB2}) of Theorem~\ref{tUB} we
need to maximize expression (\ref{eqUB}) over $d$.  The maximum
value occurs when $d=\frac{k-c-\sqrt{c(k-c)}}{k-2c}$ and is
$$ \frac{k-2\sqrt{c(k-c)}}{(k-2c)^2}\bin{n}{2}=
   \left(\frac{1}{1+2\sqrt{\frac{c}{k}
   \left(1-\frac{c}{k}\right)}}\right)\frac{1}{k}\bin{n}{2} . $$
The expression in parentheses is at most $1$. \proofend

\subsection{Proof of Corollary \ref{maincor}}

Let $H$ be a self-complementary graph with $\chib(H)=k+1$ such
that $c_{\min}$ and $c_{\max}$ are defined as in preparation for
Theorem \ref{tUB}.  That is, $c_{\min}$ is the least $c$ so that
$G$ cannot be partitioned into $c$ cliques and $k-c$ cocliques.
The quantity $c_{\max}$ is the greatest such $c$.

Because $\overline{H}=H$, $H$ can be partitioned into $c$ cliques
and $k-c$ cocliques if and only if $H$ can be partitioned into
$k-c$ cliques and $c$ cocliques.  Hence, $c_{\max}=k-c_{\min}$ and
it must be the case that $c_{\min}\leq k/2\leq c_{\max}$.  Now,
from Theorem \ref{tLB} and the first inequality of Theorem
\ref{tUB}, the result follows. \proofend



\section{Binary chromatic number}
\label{chromnum}

It is easy to see the following.
\begin{fact}  Let $G$ be a graph.
\begin{enumerate}
   \item $\chib(G)\geq \chi(G),\chi(\comp{G})$ \label{LB}
   \item $\chib(G)=\chib(\comp{G})$.
\end{enumerate}
\label{fCHIB}
\end{fact}

\noindent Recall that $K_p^q$ is the complete $p$-partite graph
with $q$ vertices in each part.

\begin{proposition}  Let $G$ be a graph.
   $$ \chib(G)\leq\chi(G)+\chi(\comp{G})-1 . $$
   This bound is tight for $G=K_p^q$.
   \label{chibUB}
\end{proposition}

\begin{proof}
Consider $c$ cliques spanning a set $A$ of vertices in $G$.
If $c<\chi(\comp{G})$
we are done since $\chi(G-A)\leq \chi(G)$. Otherwise,
$c=\chi(\comp{G})$ and it is possible to partition all vertices
into $c$ cliques. We can obtain required cocliques by considering
single vertices.

We see that $\chib(K_p^q)\geq p+q-1$ by observing that if we
require $q-1$ cliques in a partition of a vertex set of  $K_p^q$
into cliques and cocliques then $p-1$ cocliques is not enough to
partition the rest of the vertices.
\end{proof}

Next we determine the binary chromatic number of some classes of
graphs to partition the rest of the vertices.

\begin{proposition} Let $\chib(G)$ denote the binary chromatic
number of a graph $G$.
\begin{enumerate}
   \item If $n\geq 5$, then $\chib(C_n)=\lceil n/2\rceil$.
   \item If $n\geq 3$, then $\chib(P_n)=\lceil n/2\rceil$.
   \item $\chib(K_p^q)=p+q-1$. \label{kpq}
\end{enumerate}
\label{basicgraphs}
\end{proposition}

\begin{proof}
\begin{enumerate}
\item  The lower bound follows from Fact \ref{fCHIB}(\ref{LB}).
For the upper bound, we can construct the partition of a vertex
set in at most $\lceil \frac{n}{2} \rceil$ cliques and cocliques
as follows. If we need only cliques, or only cocliques, it is
clear. When we need at least one clique and at least one coclique
in that partition, take the largest coclique on $\lfloor
\frac{n}{2} \rfloor$ vertices. The leftover graph consists of
independent vertices and, if $n$ is odd, of one edge. Take this
edge (or a single vertex when $n$ is even) as a clique of our
partition. The number of leftover vertices is $\lceil \frac{n}{2}
\rceil -2$ and we are done.

\item This is quite similar to the case of $C_n$.  We leave it to
the reader.

\item This follows from Proposition~\ref{chibUB}, since
$\chi(K_p^q)=p$ and $\chi(\comp{K_p^q})=q$.
\end{enumerate}
\end{proof}

Proposition~\ref{chibLBUB} gives the bounds on the smallest binary
chromatic number among all $n$-vertex graphs.

\begin{proposition}
   If $n$ is a positive integer, then
     $$ \sqrt{n}\leq \min_{|V(G)|=n}\chib(G)\leq
   \sqrt{n}+(1+o(1))n^{0.2625}. $$
   Moreover there are infinitely many graphs for which the lower
   bound is attained.
\label{chibLBUB}
\end{proposition}

\begin{proof}
For the lower bound, we use Fact~\ref{fCHIB}(1) and the fact that
$\chi(G)\chi(\comp{G})\geq n$.  As a result, one of
$\chi(G),\chi(\comp{G})$ is larger than $\sqrt{n}$.

The lower bound is, in fact, attained by an infinite class of
graphs on $n=k^2$ vertices where $k$ is a prime. To realize this
lower bound, consider the following construction of a graph
$G=G_n=G_{k^2}$.

Let $V(G_n)$ be pairs of integers $(i,j)$, for
$i,j=1, \ldots, k$. We create $k+1$ distinct partitions of
$V(G)$ into sets of cardinalities $k$. Let the $i^{\rm th}$
partition $P_i=\{V^i_1, V^i_2, \ldots, V^i_k\}$ be defined as
follows for $i=0, \ldots, k$: $V^i_j= \{(j,1), (j+i, 2), (j+2i,
3), \ldots, (j+(k-1)i, k)\}$. Here, addition is taken modulo $k$.
Let $V^i_j$ induce a clique  if
$i<j$  and let $V^i_j$ induce a coclique if $i\geq j$.
Next we verify that $G$ is well defined.


Note that for each pair of vertices $x,y\in V(G)$, $x,y\in V_i^j$ for some $i,j$.
Moreover, if $x,y\in V^i_j$ then at most one vertex $x$ or $y$ is in $V^{i'}_{j'}$
where $i\neq i'$. Indeed, if $x,y \in V_i^j$  and $x=(x_1,x_2)$
then $y=(x_1+ li, x_2+l)$. If $x,y\in V_{i'}^{j'}$ then $y=(x_1+
l'i', x_2+l')$. Now, since $x_2+l=x_2+l'$ we have $l=l'\pmod{k}$.
Thus $x_1+li = x_1+l'i = x_1 + l'i'$, therefore $l'i=l'i'$ and
$i=i'$ if $k$ is prime.

We see that $P_i$ provides a vertex-partition of  $G$ into $i$
cocliques and $k-i$ cliques,
$0\leq i \leq k$. Therefore $\chi_B(G)\leq k=\sqrt{n}$.

For arbitrary $n$, we find the upper bound by taking the smallest
$k\geq\sqrt{n}$ such that $k$ is a prime.  Consider $G_{k^2}$ as
defined above and let $G_n$ be a subgraph  of $G_{k^2}$ induced by a
set of $n$ vertices. As we have shown, $\chib(G_{k^2})\leq k$,
which implies $\chib(G_n)\leq k$. In a paper of Baker, Hartman and
Pintz~\cite{BHP}, for $x$ at least some $x_0$, there is a prime in
the interval $\left[x-x^{0.525},x\right]$. Thus, $\chib(G)\leq
k\leq \sqrt{n}+(1+o(1))n^{0.2625}$.
\end{proof}


\section{Better bounds for small graphs}
\label{exact}

The results stated in Section~\ref{intro} are asymptotic. However,
for some graphs $H$ we are able to determine the exact value of
$\Dist(n,\Forb(H))$.

Here we shall use the fact that the extremal graphs for forbidden
induced subgraphs on three vertices as well as for induced
subgraphs on 4 vertices and 3 edges are known precisely
\cite{EFPS}.

\begin{theorem}
If $H\in \{ K_3, \comp{K_3}, K_{1,2}, \comp{K_{1,2}}\}$, then
$\Dist (n, \Forb(H))=\binom{\cn2 }{2}+ \binom{\fn2 }{2}$.
\label{threevertex}
\end{theorem}

\begin{proof}
The cases of the triangle $K_3$ and of the empty graph
$\comp{K_3}$ follow immediately from equation (\ref{turanex}).

Now, we consider the editing distance for $K_{1,2}$-free graphs.
Note that the graph which contains no induced $K_{1,2}$ is a
disjoint union of cliques.

Let $G$ be an arbitrary graph on $n$ vertices. If $G$ has minimum
degree at least $\cn2 $ then we add all missing edges to obtain a
complete graph. In this case, at most $\binom{n}{2}- \cn2
\frac{n}{2}\leq \binom{\cn2 }{2}+\binom{\fn2 }{2}$ edges were
added. Otherwise, delete all edges incident to a vertex $v$ of
degree at most $\fn2 $ and apply induction to $G\setminus v$. The
total number of additions and deletions is at most $\fn2 +
\binom{\left\lceil (n-1)/2 \right\rceil}{2}+\binom{\left\lfloor
(n-1)/2\right\rfloor }{2}\leq \binom{\cn2 }{2}+\binom{\fn2 }{2}.$
This provides an upper bound on $f$.

For the lower bound, we consider a complete bipartite graph $H$ on
$n$ vertices with almost equal parts $A$, $B$. Let $G$ be the
disjoint union of cliques $S_1, S_2, \ldots, S_k$ on the same
vertex set as $H$. Let $a_i= |A\cap V(S_i)|$ and $b_i=|B\cap
V(S_i)|$, for $i=1, \ldots, k$.  It is clear that the number of
editing operations performed on $H$ to obtain $G$ is
$$ s=\sum _{i=1}^{k}\left[\binom{a_i}{2}+ \binom{b_i}{2}
   + a_i(|B|-b_i) \right] . $$
This function is minimized when $a_i=b_i$ for all $i$, except
perhaps one $i\in\{1, \ldots, k\}$ such that $|a_i-b_i|=1$. Now,
$s \geq  n^2/4-n/2$ for even $n$ and $s \geq (n-1)^2/4 - (n-1)/2 +
(n-1)/2$ for odd $n$, and the result follows.
\end{proof}

Let $\cal Q$ be the set of graphs on $n$ vertices with no induced
subgraphs on $4$ vertices and $3$ edges. In \cite{EF} it was shown
that any graph in $\cal{Q}$ or its complement is a disjoint union
of $4$-cycles and trees on at most $3$ vertices.  Note that
$G\in{\cal Q}$ if and only if $\overline{G}\in{\cal Q}$.

\begin{theorem}
$\lfloor(n^2-5n)/4\rfloor \leq  \Dist(n, {\cal Q})\leq (n^2- n)/4
$. \label{4n3e}
\end{theorem}
\begin{proof}
Let $G$ be a graph on $n$ vertices.  Since $E_n, K_n\in \cal{Q}$,
it is sufficient either to add all missing edges or to delete all
edges to obtain a graph from $\cal Q$.  Thus, the upper bound
follows.

For the lower bound consider a graph $G$ with
$\lfloor(n^2-n)/4\rfloor$ edges.  Assume first that the minimum
number of edit operations results in a graph whose components are
either $4$-cycles or trees on at most $3$ vertices. As a result,
the total number of edges within these components is at most $n$.
Therefore, at least $|E(G)|-n$ edges of $G$ had to be deleted.

The result is similar if the minimum number of edit operations
results in a graph such that its complement has components that
are either $4$-cycles or trees on at most $3$ vertices.  So, at
least $|E(\overline{G})|-n=\bin{n}{2}-|E(G)|-n$ edges had to be
added to $G$.

As a result, the number of edit operations is at least
$\lfloor(n^2-n)/4\rfloor-n$.
\end{proof}

\noindent
{\bf Conclusions}

The editing problem of graphs we consider in this paper can be
reformulated in terms of complete edge-colored graphs, where the
edges of the graphs correspond to edges of one color, say red, and
the edges of the complement correspond to edges of another color,
say blue. Our editing operations are equivalent to changing the
color of some edges from red to blue or from blue to red.

It is natural to consider more than two colors. Specifically for
any two colorings of $E(K_n)$ in colors from $\{1, \ldots,
\gamma\}$, we define the distance to be the smallest number of
edge-recolorings to obtain one coloring from the other. Our
results for classes of graphs with forbidden induced subgraphs can
be generalized for classes of multicolored graphs with forbidden
color patterns.  When considered on bipartite graphs, the
multicolored graph editing problem is equivalent to the problem of
editing a matrix so that fixed patterns on submatrices do not
occur~\cite{AM}. \\

\noindent
{\bf Acknowledgements}

We thank Ling Long for alerting us to the result of Baker, Hartman
and Pintz~\cite{BHP}.


\section{Appendix: Random graph}

Let $G(n,p)$ be a
graph in which each edge from $K_n$ is chosen independently with
probability $p$ (see \cite{Bbook,JLR}).
Lemma \ref{cRANDOM} follows immediately from the following:

\begin{lemma}
   Fix a constant $\epsilon>0$ and positive integer $m$.
   Let $\ell$ have the property that $m\leq\ell\leq M(\epsilon,m)$
   where $M(\epsilon,m)$ is the constant given in the Regularity
   Lemma (Lemma \ref{reglem}).  Let $G=G(n,1/2)$,
   $f(n)=\omega(n^{-1/2})$ and $P$ be the probability that for
   each $(m,\epsilon,\ell)$-equipartition of the vertices of $G$,
   all pairs of clusters $(V_i,V_j)$, $1\leq i<j\leq \ell$,
   have density in the interval $(1/2-f(n), 1/2+f(n))$. Then
   $P$ approaches $1$ as $n$ goes to infinity.
   \label{lRANDOM}
\end{lemma}

\begin{proof}
We just want to compute the probability that {\bf all} pairs of
disjoint sets, each of size at least $\epsilon' n$ (where
$\epsilon'=\frac{1-\epsilon}{M(\epsilon)-1}$), have density in the
interval $\left(1/2-f,1/2+f\right)$, for any
$f=f(n)=\omega(n^{-1/2})$.

\begin{eqnarray}
   \lefteqn{\Pr\left\{\bigvee_{\scriptsize \begin{array}{l}
                               S,T\subseteq V(G) \\
                               S\cap T=\emptyset \\
                               |S|,|T|\geq \epsilon' n
                               \end{array}}
                      (d(S,T)\not\in (1/2-f,
                                   1/2+f))\right\}}
            \nonumber \\
   & \leq & 2^n2^n 2\Pr\left\{d(S,T)<1/2-f\right\}
            \nonumber \\
   & \leq & 2\cdot 4^n\exp\left(-2(f|S||T|)^2/(|S||T|)\right)
            \label{CHERNOFF} \\
   & \leq & 2\cdot 4^n\exp\left(-2f^2|S||T|\right) \nonumber \\
   & \leq & 2\cdot 4^n\exp\left(-2(\epsilon')^2f^2n^2\right)
            \rightarrow 0 \nonumber
\end{eqnarray}

Chernoff's bound (see~\cite{JLR}) is used to achieve inequality
(\ref{CHERNOFF}).
\end{proof}


\begin{thebibliography}{99}

\bibitem{AM} M. Axenovich and R. Martin, Avoiding patterns in
matrices via small number of changes, submitted for publication.

\bibitem{BHP} R.C. Baker, G. Hartman and J. Pintz, The difference
between consecutive primes, II. {\it Proc. London Math. Soc.} {\bf
83} (2001), no. 3, 532--562.

\bibitem{B} B. Bollob\'{a}s, Hereditary properties of graphs:
asymptotic enumeration, global structure, and colouring.
Proceedings of the international Congress of Mathematicians, Vol.
III (Berlin, 1998).  {\it Doc. Math.} 1998, Extra Vol. III,
333-342 (electronic).

\bibitem{Bbook} B. Bollob\'{a}s, {\it Random Graphs}.  Second edition.
Cambridge Studies in Advanced Mathematics, 73.  {\it Cambridge
University Press, Cambridge}, 2001.

\bibitem{BT} B. Bollob\'{a}s and A. Thomason, Hereditary and
monotone properties of graphs.  {\it The Mathematics of Paul
Erd\H{o}s, II} 70--78, Algorithms Combin., 14, Springer, Berlin,
1997.

\bibitem{CEFS} D. Chen, O. Eulenstein, D. Fern{\'a}ndez-Baca and
M. Sanderson, Supertrees by flipping, preprint.

\bibitem{DF} R.G. Downey and M.R. Fellows, {\it Parameterized
complexity}.  Springer, New York, 1999.

\bibitem{EFPS} P. Erd\H{o}s, R. Faudree, J. Pach and J. Spencer,
How to make a graph bipartite. {\it J. Combin. Theory Ser. B} {\bf
45} (1988), no. 1, 86--98.

\bibitem{EF} P. Erd\H{o}s, Z. F\"uredi, B. L. Rothschild, and
V. T. S\'os, Induced subgraphs of given sizes.  {\it Paul
Erd{\H{o}}s memorial collection.} Discrete Math. {\bf 200} (1999),
no. 1-3, 61--77.

\bibitem{EGR} P. Erd\H{o}s, A. Gy\'arf\'as and M. Ruszink\'o, How to
decrease the diameter of triangle-free graphs.  {\it
Combinatorica} {\bf 18} (1998), no. 4, 493--501.

\bibitem{EGS} P. Erd\H{o}s, E. Gy\H{o}ri and M. Simonovits, How many
edges should be deleted to make a triangle-free graph bipartite?
{\it Sets, graphs and numbers} (Budapest, 1991), 239--263, Colloq.
Math. Soc. J{\'a}nos Bolyai, {\bf 60}, North-Holland, Amsterdam,
1992.

\bibitem{ES2} P. Erd\H{o}s and M. Simonovits, A limit theorem in
graph theory.  {\it Studia Sci. Math. Hungar.} {\bf 1}, (1966),
51--57.

\bibitem{ES} P. Erd\H{o}s and A.H. Stone, On the structure of linear
graphs. {\it Bull. Amer. Math. Soc.} {\bf 52}, (1946), 1087--1091.

\bibitem{F} Z. F{\"{u}}redi, Tur{\'a}n type problems.  {\it
Surveys in combinatorics,} 253--300, London Math. Soc. Lecture
Note Ser., {\bf 166}, Cambridge Univ. Press, Cambridge, 1991.

\bibitem{JLR} S. Janson, T. {\L}uczak and A. Ruci{\'{n}}ski, {\it
Random Graphs}.  Wiley-Interscience Series in Discrete Mathematics
and Optimization.  Wiley-Interscience, New York, 2000.

\bibitem{KSSSz} J. Koml\'os, A. Shokoufandeh, M. Simonovits and
E. Szemer\'edi, The regularity lemma and its applications in graph
theory.  {\it Theoretical aspects of computer science (Tehran,
2000),} 84--112, Lecture Notes in Comput. Sci., {\bf 2292},
Springer, Berlin, 2002.

\bibitem{KS} J. Koml\'os and M. Simonovits, Szemer\'edi's regularity
lemma and its applications in graph theory. {\it Combinatorics,
Paul Erd\H{o}s is eighty, Vol. 2 (Keszthely, 1993),} 295--352,
Bolyai Soc. Math. Stud. 2, J{\'a}nos Bolyai Math. Soc., Budapest,
1996.

\bibitem{LS} L.M. Lesniak and J.H. Straight, The
cochromatic number of a graph.  {\it Ars Combin.} {\bf 3} (1977),
39--45.


\bibitem{PS} H.J. Pr\"omel and A. Steger, Excluding induced
subgraphs. III A general asymptotic.  {\it Random Structures
Algorithms} {\bf 3} (1992), no. 1, 19--31.

\bibitem{S} M. Simonovits, Extremal graph theory,
{\it Selected topics in graph theory}, {\bf 2}, 161--200, Academic
Press, London, 1983.

\bibitem{SAN1998} G. Stephanopoulos, A. Aristidou and J. Nielsen,
{\it Metabolic engineering: principles and methodologies}.
Academic Press, San Diego, 1998.

\bibitem{Sz} E. Szemer\'edi, Regular partitions of graphs.  In
{\it Probl{\`e}mes Combinatoires et Th{\'e}orie des Graphes},
399--401, Colloq. Internat. CNRS, Univ. Orsay, Paris, 1978.

\bibitem{T} P. Tur\'{a}n, Eine Extremalaufgabe aus der
Graphentheorie.  {\it Mat. Fiz. Lapok} {\bf 48}, (1941), 436--452.

\bibitem{W} D. West, {\it Introduction to graph theory}. Second edition.
{\it Prentice Hall, Inc., Upper Saddle River, NJ, 1996}, p. 588.

\end{thebibliography}
\end{document}